\documentclass[11pt]{article}
\usepackage{amsmath,amssymb,amsfonts}
\usepackage{amsthm}
\usepackage{cite}
\newtheorem{theorem}{Theorem}

\newtheorem{lemma}{Lemma}

\title{\textbf{Central Limit Theorem for a Class of SPDEs} }
\date{}
\author{Parisa Fatheddin \\
\footnotesize University of Tennessee, Knoxville, TN 37996, USA
}
\begin{document}
\newtheorem{example}[theorem]{Example}
\newtheorem{cor}[theorem]{Corollary}
\newtheorem{notation}[theorem]{Notation}
\newtheorem{notations}[theorem]{Notations}
\newtheorem{claim}[theorem]{Claim}
\newtheorem{mtheorem}[theorem]{Meta-Theorem}
\newtheorem{prop}[theorem]{Proposition}
\newtheorem{rem}[theorem]{Remark}
\newtheorem{conj}[theorem]{Conjecture}
\newtheorem{rems}[theorem]{Remarks}
\maketitle

\begin{abstract}
Here we establish the central limit theorem for a class of stochastic partial differential equations (SPDEs) and as an application derive this theorem for two widely studied population models known as super-Brownian motion and Fleming-Viot process.
\end{abstract}

\noindent {\sc Mathematics Subject Classification (2010):} Primary 60F05;
Secondary: 60H15, 60J68.

\noindent {\sc Key words:} Central limit theorem, stochastic
partial differential
 equation, Fleming-Viot process, super-Brownian motion.

\section{Introduction}
Two commonly studied population models are super-Brownian motion (SBM) and Fleming-Viot Process (FVP). These are measure-valued Markov processes and can be represented as SPDEs. We use these representations to formulate a general class of SPDE and investigate the central limit theorem for this class and have the two population models as special cases. These models are formed as the scaled limit of their discrete particle systems.

SBM is the continuous version of the Branching Brownian motion, the oldest and best known branching process and individuals are assumed to reproduce following a Galton-Watson process. In this model, the population evolves as a ``cloud'' through time with each individual  assumed to move according to a Brownian motion and leave behind a random number of offsprings upon death. On the other hand, FVP is the continuous approximation of step-wise mutation process, in which each individual has a ``type'' (usually genetic type) given by an element $x$ in some set $E$. In this model, we are interested in the distributions of the types in the whole population making FVP a probability measure-valued process. Mutation is the term referring to a change in genetic type. In FVP the number of individuals is assumed to be fixed throughout time;that is, in the place of an individual's death an offspring is born. For more information and background on SBM and FVP and their formulation as the continuous approximation of discrete particle systems, we refer the reader to \cite{Eth}.

On the topic of central limit theorem (CLT), developments have been made on SBM and on various processes related to it. Li \cite{Li} considered the critical continuous SBM and proved the CLT in all dimensions, $d\geq 1$, and also derived the CLT for its weighted occupation time process in $d\geq 3$. Schied showed the tightness and weak convergence of the finite-dimensional distributions of SBM to those of a Wiener process as part of the proof of the moderate deviation principle for SBM in \cite{Schied} for $d\geq 1$. CLT for SBM was then concluded for all dimensions. In \cite{Lee}, Lee and Remillard also used their large and moderate deviation results to derive the CLT for SBM in dimension $d=3$. In addition, some authors have studied SBM with super-Brownian immigration (SBMSBI). Hong and Li \cite{Li1} proved the CLT for SBMSBI for dimensions $d\geq 3$ and later Hong \cite{HongQ} showed its CLT under the quenched probability law for the same dimensions. With Zeitouni, Hong also succeeded in achieving the quenched CLT for SBMSBI for $d\geq 4$ in \cite{HongZ}. In addition, Hong showed the CLT for the occupation time process of SBMSBI for $d\geq 3$ in \cite{Hong}. The CLT for SBM with other types of immigration have also been considered. See for example \cite{Li2, ZhangW, HongC,HongW}. To the knowledge of the author, CLT has not been previously shown for FVP.

All authors mentioned above proved the CLT by using the Laplace transform of the process under study. Except Schied \cite{Schied}, Lee and Remillard \cite{Lee} and Zhang \cite{ZhangW}, they applied the method offered by Iscoe \cite{Iscoe} to achieve the CLT. To be more precise, Iscoe's method consists of finding the limit of the Laplace functional of the centered process and applying the Bochner-Minlos theorem. That is, to achieve the CLT for a process, $\{X_{t}\}$, first the centered functional process is formed by
\begin{equation*}
\left<Z_{t},f\right>= a_{d}^{-1}(t)\left(\left<X_{t},f\right>-\mathbb{E}\left<X_{t},f\right>\right)
\end{equation*}
for some norming constant, $a_{d}(t)$, and for $f\in \mathcal{S}(\mathbb{R})$, the Schwartz space. Then the weak convergence of $Z_{t}$ to a centered Gaussian process, $Z_{\infty}$, is obtained by considering the Laplace functional of the centered functional process,
\begin{equation*}
\mathbb{E}\exp\left(-\left<Z_{t},f\right>\right)= \exp(W_{t})
\end{equation*}
where the limit of $W_{t}$ is the covariance of the Gaussian process $Z_{\infty}$.

 Here instead of the Laplace transform of SBM, we use another characterization of this population model given by Xiong \cite{Xio}. In \cite{Xio}, by studying SBM as a ``distribution'' function-valued process, an SPDE was formed to define SBM. A similar SPDE was also derived for FVP. By observing the similarities between the two SPDEs we formulate a general SPDE and derive the CLT in $d=1$ for this SPDE and as an application establish the CLT for the two population models. We note that since the formulation of the general SPDE and the two population models given in \cite{Xio} is in $d=1$ only, our results are limited to this dimension. Extending the result of \cite{Xio} and also the result presented here to higher dimensions require further investigation.

We begin by some background and notations in Section 2. We then prove the CLT for the general SPDE in Section 3 by first showing its tightness in our space notified in Section 2 and afterwards proving that the limiting process has a unique solution and is Gaussian. Section 4 contains the CLT for the two population models, SBM and FVP.

\section{Notations and Main Results}
Suppose $(\Omega, \mathcal{F}, P)$ is a probability space and $\{\mathcal{F}_{t}\}$ is a family of non-decreasing right continuous sub-$\sigma$-fields of $\mathcal{F}$ such that $\mathcal{F}_{0}$ contains all $P$-null subsets of $\Omega$. We denote $\mathcal{C}_{b}(\mathbb{R})$ to be the space of continuous bounded functions on $\mathbb{R}$ and $\mathcal{C}_{c}(\mathbb{R})$ to be composed of continuous functions in $\mathbb{R}$ with compact support. Let $K$ be a constant that may change values at different lines. Let
\begin{equation*}
\mathcal{C}_{p}(\mathbb{R}^{d}):= \left\{f\in \mathcal{C}(\mathbb{R}^{d}): \sup\frac{|f(x)|}{\phi_{p}(x)} <\infty \mbox{ for } p>d, \phi_{p}(x):=(1+|x|^{2})^{-\frac{p}{2}}\right\}
\end{equation*}
Since SBM is a measure-valued process, we denote it by $\mu_{t}^{\epsilon}$ with branching rate $\epsilon$. There are two common ways to define SBM, $\mu_{t}^{\epsilon}$. One is by its Laplace transform given by,
\begin{equation*}
\mathbb{E}_{\mu_{0}^{\epsilon}} \exp(-<\mu_{t}^{\epsilon},f>)= \exp(-<\mu_{0}^{\epsilon}, v(t,\cdot)>)
\end{equation*}
where $v(\cdot, \cdot)$ is the unique mild solution of the evolution equation:
\begin{equation*}
\left\{\begin{array}{l} \dot{v}(t,x)= \frac{1}{2}\Delta v(t,x) - v^{2}(t,x)  \\
v(0,x)=f(x)
\end{array}\right.
\end{equation*}
for $f\in \mathcal{C}_{p}^{+}(\mathbb{R}^{d})$ and the other is as the unique solution to a martingale problem: for all $f\in \mathcal{C}_{b}^{2}(\mathbb{R})$
\begin{equation*}
M_{t}(f):= \left<\mu_{t}^{\epsilon},f\right>-\left<\mu_{0}^{\epsilon},f\right>-\int_{0}^{t}\left<\mu_{s}^{\epsilon},\frac{1}{2}\Delta f\right> ds
\end{equation*}
is a square-integrable martingale with quadratic variation,
\begin{equation*}
<M(f)>_{t}= \epsilon \int_{0}^{t} <\mu_{s}^{\epsilon},f^{2}>ds
\end{equation*}
Similarly, let $\mu_{t}^{\epsilon}$ denote FVP with mutation rate $\epsilon$. Here $\{\mu_{t}^{\epsilon}\}$ is a family of probability measures and there are also two usual ways of defining this process. One is as a Markov process with generator,
\begin{eqnarray*}
\mathcal{L}^{\epsilon}F(\mu_{t}^{\epsilon})&=& f'(<\mu_{t}^{\epsilon},\phi>)<\mu_{t}^{\epsilon},A\phi> \\
&&\hspace{.4cm}+ \frac{\epsilon}{2}\int\int f''(<\mu_{t}^{\epsilon},\phi>)\phi(x)\phi(y)Q(\mu_{t};dx,dy)
\end{eqnarray*}
having domain,
\begin{equation*}
\mathcal{D}= \{F: F(\mu_{t}^{\epsilon})=f(<\mu_{t}^{\epsilon},\phi>),f\in \mathcal{C}_{b}^{\infty}(\mathbb{R}), \phi \in D(A),\mu \in M_{1}(E)\}
\end{equation*}
where $\mathcal{C}_{b}^{\infty}(\mathbb{R})$ is the set of all bounded, infinitely differentiable functions on $\mathbb{R}$, $M_{1}(E)$ is the space of all probability measures on $E$ endowed with the usual weak topology. Furthermore, $D(A)$ denotes the domain of $A$, where $A$ is the generator of a Markov process on the set $E=[0,1]$. In the context of population models, $E$ represents the genetic type space of the population and $A$ is referred to as the mutation operator. Moreover,
\begin{equation*}
Q(\mu_{t}^{\epsilon};dx,dy)= \mu_{t}^{\epsilon}(dx)\delta_{x}(dy)-\mu_{t}^{\epsilon}(dx)\mu_{t}^{\epsilon}(dy)
\end{equation*}
where $\delta_{x}$ denotes the Dirac measure at $x$. For more information on this characterization of FVP see \cite{DF1} and \cite{FenX}.

 The second way to define FVP is as a unique solution to a martingale problem: for $f\in \mathcal{C}_{c}^{2}(\mathbb{R})$,
\begin{equation*}
M_{t}(f)= <\mu_{t}^{\epsilon},f>-<\mu_{0}^{\epsilon},f>-\int_{0}^{t} <\mu_{s}^{\epsilon},\frac{1}{2}\Delta f>ds
\end{equation*}
is a continuous square-integrable martingale with quadratic variation,
\begin{equation*}
<M_{t}(f)>= \epsilon \int_{0}^{t}\left(<\mu_{s}^{\epsilon},f^{2}>-<\mu_{s}^{\epsilon},f>^{2}\right)ds
\end{equation*}
Recently, another formulation of SBM and FVP was given in \cite{Xio} by considering their ``distribution'' function-valued process. More precisely, by considering $u_{t}^{\epsilon}(y)= \int_{0}^{y}\mu_{t}^{\epsilon}(dx)$ for all $y\in \mathbb{R}$, SBM was characterized in \cite{Xio} by the following SPDE:

\begin{equation}\label{SBM}
u_{t}^{\epsilon}(y)=F(y)+ \int_{0}^{t}\int_{0}^{u_{s}^{\epsilon}(y)} W(dsda)+  \int_{0}^{t} \frac{1}{2} \Delta u_{s}^{\epsilon}(y)ds
\end{equation}
where $F(y)= \int_{0}^{y}\mu_{0}^{\epsilon}(dx)$ and $W$ is a white noise random measure on $\mathbb{R}^{+}\times \mathbb{R}$ with intensity measure $dsda$.

Also by considering $u_{t}^{\epsilon}(y)=\int_{-\infty}^{y}\mu_{t}^{\epsilon}(dx)$ for all $y \in \mathbb{R}$, FVP was given in \cite{Xio} by the SPDE,
\begin{equation}\label{FVP}
u_{t}^{\epsilon}(y)= F(y) + \int_{0}^{t}\int_{0}^{1} \left(1_{a\leq u_{s}^{\epsilon}(y)}-u_{s}^{\epsilon}(y)\right)W(dsda) + \int_{0}^{t} \frac{1}{2}\Delta u_{s}^{\epsilon}(y)ds
\end{equation}
We denote both processes as $\mu_{t}^{\epsilon}$ where based on the context it will be clear which process is being referred to. By noticing the similarities between the two SPDE formulations given above, we form a general SPDE, with the two models as special classes, as follows,
\begin{equation}\label{SPDE}
u_{t}^{\epsilon}(y)=F(y) + \sqrt{\epsilon} \int_{0}^{t}\int_{U}G(a,y,u_{s}^{\epsilon}(y)) W(dsda) + \int_{0}^{t} \frac{1}{2}\Delta u_{s}^{\epsilon}(y)ds
\end{equation}
where $G: U \times \mathbb{R}^{2} \rightarrow \mathbb{R}$, $F$ is a function on $\mathbb{R}$ and for $u_{1},u_{2},u,y\in \mathbb{R}$,
\begin{eqnarray}
\int_{U}\left|G(a,y,u_{1})-G(a,y,u_{2})\right|^{2} \lambda(da) &\leq& K|u_{1}-u_{2}| \label{con1}\\
\int_{U}|G(a,y,u)|^{2}\lambda(da) &\leq& K(1+|u|^{2})\label{con2}.
\end{eqnarray}

Let $\mathcal{S}(\mathbb{R})$ be the Schwartz space of rapidly decreasing functions defined as
\begin{equation*}
\mathcal{S}(\mathbb{R})= \left\{\phi \in \mathcal{C}^{\infty}(\mathbb{R}): \|\phi\|_{\alpha, \beta} < \infty, \forall \alpha, \beta \in \mathbb{N}\cup \{0\}\right\}
\end{equation*}
where
\begin{equation*}
\|\phi\|_{\alpha,\beta} = \sup_{x\in \mathbb{R}} \left|x^{\alpha}\phi^{(\beta)}(x)\right|
\end{equation*}
with its dual, $\mathcal{S}'(\mathbb{R})$, known as the space of tempered distributions. To investigate the CLT for the general SPDE, we consider the $\mathcal{S}'(\mathbb{R})$-valued centered process:
\begin{equation}\label{form}
Z_{t}^{\epsilon}= \frac{1}{\sqrt{\epsilon}}( u^{\epsilon}_{t}-u^{0}_{t})
\end{equation}
Namely, we study the process:
\begin{equation}\label{eq1}
<Z_{t}^{\epsilon},f> :=\int_{0}^{t} \int_{U} \int_{\mathbb{R}} G(a,y,u_{s}^{\epsilon}(y)) f(y)dyW(dads) + \frac{1}{2} \int_{0}^{t} <Z_{s}^{\epsilon},f''>ds
\end{equation}
for $f\in \mathcal{S}(\mathbb{R})$.

\begin{theorem}
The centered process, $\{Z_{t}^{\epsilon}\}$, is tight in $\mathcal{C}\left([0,1]; \mathcal{S}'(\mathbb{R})\right)$.
\end{theorem}
We use the above theorem to obtain the following results on CLT.

\begin{theorem}
The general SPDE, $\{u_{t}^{\epsilon}\}$, satisfies the CLT in space $\mathcal{C}\left([0,1], \mathcal{S}'(\mathbb{R})\right)$, where $\{Z_{t}^{\epsilon}\}$ converges in distribution as $\epsilon$ tends to zero to a Gaussian process, $\{Z_{t}^{0}\}$ with zero mean and covariance,
\begin{eqnarray}\label{co1}
&&Cov\left(\left<Z_{t}^{0},f\right>,\left<Z_{t}^{0},g\right>\right)\\
&=& \int_{0}^{t}\int_{U}\int_{\mathbb{R}} G(a,y,u_{s}^{0}(y))f(y)dy\int_{\mathbb{R}} G(a,x,u_{s}^{0}(x))g(x)dx\lambda(da)ds\nonumber
\end{eqnarray}
for $f,g \in \mathcal{S}(\mathbb{R})$.
\end{theorem}
For the next two theorems, let $\mu_{t}^{\epsilon}$ denote SBM and FVP and consider the centered process,
\begin{equation*}
\tilde{Z}_{t}^{\epsilon} = \frac{1}{\sqrt{\epsilon}}\left(\mu_{t}^{\epsilon}-\mu_{t}^{0}\right)
\end{equation*}

\begin{theorem}
 SBM satisfies the CLT in space $\mathcal{C}\left([0,1], \mathcal{S}'(\mathbb{R})\right)$, where $\{\tilde{Z}_{t}^{\epsilon}\}$ converges in distribution as $\epsilon$ tends to zero to a Gaussian process, $\{\tilde{Z}_{t}^{0}\}$ with zero mean and covariance,
\begin{equation}\label{co2}
Cov\left(\left<\tilde{Z}_{t}^{\epsilon},f\right>,\left<\tilde{Z}_{t}^{\epsilon},g\right>\right)= \int_{0}^{t}
\left<\mu_{s}^{0},fg\right>ds
\end{equation}
for $f,g \in \mathcal{S}(\mathbb{R})$.
\end{theorem}

\begin{theorem}
  FVP satisfies the CLT in space $\mathcal{C}\left([0,1], \mathcal{S}'(\mathbb{R})\right)$, where $\{\tilde{Z}_{t}^{\epsilon}\}$ converges in distribution as $\epsilon$ tends to zero to a Gaussian process, $\{\tilde{Z}_{t}^{0}\}$ with zero mean and covariance,
\begin{eqnarray}\label{co3}
&&Cov\left(\left<\tilde{Z}_{t}^{\epsilon},f\right>, \left<\tilde{Z}_{t}^{\epsilon},g\right>\right)\\
&=&\int_{0}^{t}\int_{0}^{1}f(y)g(y)\mu_{s}^{0}(dy)ds - \int_{0}^{t}\int_{0}^{1}\left<\mu_{s}^{0},f\right>g(y)\mu_{s}^{0}(dy)ds\nonumber\\
&& \hspace{.3cm} - \int_{0}^{t}\int_{0}^{1}f(y)\left<\mu_{s}^{0},g\right>\mu_{s}^{0}(dy)ds + \int_{0}^{t}\left<\mu_{s}^{0},f\right>\left<\mu_{s}^{0},g\right>ds\nonumber
\end{eqnarray}
for $f,g \in \mathcal{S}(\mathbb{R})$.
\end{theorem}

\section{CLT for the General SPDE}
We begin by proving Theorem 1. Since strong uniqueness of solutions to general SPDE $\{u_{t}^{\epsilon}\}$ was obtained in \cite{Xio}, then there exists a unique solution to $Z_{t}^{\epsilon}$;consequently, we have the uniqueness of solutions to our process of study, $\left<Z_{t}^{\epsilon},f\right>$. Thus, we use its mild solution instead, given by
\begin{equation}\label{mild}
\left<Z_{t}^{\epsilon},f\right> = \int_{0}^{t}\int_{U}\int_{\mathbb{R}} P_{t-s}G(a,y,u_{s}^{\epsilon}(y))f(y)dyW(dads)
\end{equation}
where $P_{t-s}$ is the Brownian semigroup defined as $P_{t}f(x)=\int_{\mathbb{R}} p_{t}(x-y)f(y)dy$ with $p_{t}(x-y)=\frac{1}{\sqrt{2\pi t}} e^{-\frac{|x-y|^{2}}{2t}}$ being the heat kernel. We show that $\left<Z_{t}^{\epsilon},f\right>$ is tight in $\mathcal{C}\left([0,1]; \mathcal{S}'( \mathbb{R})\right)$ by applying a classic result given below.

\begin{theorem}[Theorem $12.3$ in \cite{Billi}]\label{T12.3}
The sequence $\{X_{n}\}$ is tight in $\mathcal{C}\left([0,1];\mathbb{R}\right)$, if it satisfies these two conditions:\\
$(i)$ The sequence $\{X_{n}(0)\}$ is tight\\
$(ii)$ There exist constants $\gamma \geq 0$ and $\alpha >1$ and a nondecreasing, continuous function $F$ on $[0,1]$ such that \begin{equation}\label{12.3}
P\left(\left|X_{n}(t_{2})-X_{n}(t_{1})\right|\geq \lambda\right) \leq \frac{1}{\lambda^{\gamma}}\left|F(t_{2})-F(t_{1})\right|^{\alpha}
\end{equation}
holds for all $t_{1},t_{2}$ and $n$ and all positive $\lambda$.
\end{theorem}
As stated in \cite{Billi}, the moment condition,
\begin{equation}\label{moment}
\mathbb{E}\left(\left|X_{n}(t_{2})-X_{n}(t_{1})\right|^{\gamma}\right)\leq \left|F(t_{2})-F(t_{1})\right|^{\alpha}
\end{equation}
implies (\ref{12.3}). Note that in our case, $\{X_{n}(0)\}=0$  for all $n$, therefore it is sufficient to prove condition two in above theorem by checking that $\left<Z_{t}^{\epsilon},f\right>$ satisfies the moment condition given by (\ref{moment}). To this end, we use the following lemma, the proof of which is very similar to the one given for Lemma 2.3 in \cite{Xio}.
\begin{lemma}\label{Mlemma}
For any $n\geq 2$,
\begin{equation*}
\mathbb{E}\left(\sup_{\epsilon > 0}\sup_{0\leq s\leq 1} \int_{\mathbb{R}} |u_{s}^{\epsilon}(x)|^{2}e^{-2|x|}dx\right)^{n}<\infty
\end{equation*}
\end{lemma}
\begin{flushright}
  $\Box$
\end{flushright}

Given any $t_{1},t_{2}\in [0,1]$ without loss of generality we assume $t_{1}<t_{2}$ and with the help of Burkholder-Davis-Gundy inequality and condition (\ref{con2}) we obtain,
\begin{eqnarray}\label{Zt}
&&\mathbb{E}\left|\left<Z_{t_{2}}^{\epsilon},f\right>-\left<Z_{t_{1}}^{\epsilon},f\right>\right|^{4}\\
&=& \mathbb{E}\left|\int_{t_{1}}^{t_{2}}\int_{U}\int_{\mathbb{R}}P_{t-s} G(a,y,u_{s}^{\epsilon}(y))f(y)dyW(dads)\right|^{4}\nonumber\\
&\leq& \mathbb{E}\left(\int_{t_{1}}^{t_{2}}\int_{U}\left(\int_{\mathbb{R}} P_{t-s}G(a,y,u_{s}^{\epsilon}(y))f(y)dy\right)^{2}\lambda(da)ds\right)^{2}\nonumber\\
&\leq& \mathbb{E}\left(\int_{t_{1}}^{t_{2}}\int_{U}\int_{\mathbb{R}}G(a,y,u_{s}^{\epsilon}(y))^{2}e^{-2|\lfloor y\rfloor|}dy
\int_{\mathbb{R}} \left(P_{t-s}f(r)\right)^{2}e^{2|\lfloor r\rfloor|}dr\lambda(da)ds\right)^{2}\nonumber\\
&\leq& \mathbb{E}\left(\int_{t_{1}}^{t_{2}}\int_{\mathbb{R}} K\left(1+\left|u_{s}^{\epsilon}(y)\right|^{2}\right)e^{-2|\lfloor y\rfloor|}dy\int_{\mathbb{R}}\left(P_{t-s}f(r)\right)^{2}
e^{2|\lfloor r\rfloor|}drds\right)^{2}\nonumber
\end{eqnarray}
where $\lfloor y\rfloor$ is the greatest integer less than or equal to $y$. Using the fact that for any $f\in \mathcal{S}(\mathbb{R}^{d})$,
\begin{equation*}
\|P_{s}f\| \leq K\left(1\wedge s^{-d/2}\right)
\end{equation*}
we obtain,
\begin{eqnarray*}
&&\int_{\mathbb{R}} \left(P_{t-s}f(y)\right)^{2}e^{2|\lfloor y\rfloor|}dy\\
&\leq& K\left(1\wedge (t-s)^{-1/2}\right) \int_{\mathbb{R}} P_{t-s}f(y)e^{2|\lfloor y\rfloor|}dy\\
&\leq& K\left(1\wedge (t-s)^{-1/2}\right)\int_{\mathbb{R}}\int_{\mathbb{R}} e^{-\frac{|x|^{2}}{2(t-s)}}\left(e^{\lceil \frac{2|xy|}{2(t-s)}\rceil}f(x)\right)e^{-\frac{|y|^{2}}{2(t-s)}+2|\lfloor y\rfloor|}dxdy\\
&\leq& K
\end{eqnarray*}
where $\lceil y\rceil$ is the least integer greater than or equal to $y$. Hence by Lemma \ref{Mlemma},
\begin{eqnarray*}
\mathbb{E}\left|\left<Z_{t_{1}}^{\epsilon},f\right>-\left<Z_{t_{2}}^{\epsilon},f\right>\right|^{4}\leq K\left|t_{1}-t_{2}\right|^{2}
\end{eqnarray*}
Therefore, $\{\left<Z_{t}^{\epsilon},f\right>\}$ is tight in $\mathcal{C}\left([0,1]; \mathbb{R}\right)$. Based on Theorem 3.1 of \cite{Mit} we conclude the tightness of $\left\{\left<Z_{t}^{\epsilon},f\right>\right\}$ in $\mathcal{C}\left([0,1]; \mathcal{S}'(\mathbb{R})\right)$.

Now we prove the convergence of the centered process of the general SPDE to a Gaussian process to obtain the CLT for the general SPDE. We can achieve the $L^{2}$-convergence of $\{Z_{t}^{\epsilon}\}$ as $\epsilon \rightarrow 0$ using the Burkholder-Davis-Gundy inequality and condition (\ref{con1}) of $G(a,.,u_{s}^{\epsilon}(.))$ as follows.
\begin{eqnarray}\label{unique}
&&\mathbb{E}\left|Z_{t}^{\epsilon}(y)-Z_{t}^{0}(y)\right|^{2}\\
&=& \mathbb{E}\left|\int_{0}^{t}\int_{U}P_{t-s}\left(G(a,y,u_{s}^{\epsilon}(y))-G(a,y,u_{s}^{0}(y))\right)W(dads)\right|^{2}\nonumber\\
&\leq& \mathbb{E}\left|\int_{0}^{t}\int_{U}\left(P_{t-s}\left(G(a,y,u_{s}^{\epsilon}(y))-G(a,y,u_{s}^{0}(y))\right)\right)^{2}
\lambda(da)ds\right|\nonumber\\
&\leq& K\mathbb{E}\left|\int_{0}^{t}P_{2t-2s}\left|u_{s}^{\epsilon}(y)-u_{s}^{0}(y)\right|
ds\right|\nonumber
\end{eqnarray}
Note that
\begin{eqnarray*}
&&\left|\int_{0}^{t} P_{2t-2s}\left|u_{s}^{\epsilon}(y)-u_{s}^{0}(y)\right|ds\right|\\
&=& \left|\int_{0}^{t}\int_{\mathbb{R}}p_{2t-2s}(x-y)\left|u_{s}^{\epsilon}(x)-u_{s}^{0}(y)\right|dxds\right|\\
&\leq& \left|\int_{0}^{t}\left(\int_{\mathbb{R}}p_{2t-2s}^{2}(x-y)e^{2|x|}dx\right)^{1/2}\left(\int_{\mathbb{R}} \left|u_{s}^{\epsilon}(x) - u_{s}^{0}(x)\right|^{2}e^{-2|x|}dx\right)^{1/2}ds\right|\\
&\leq& K
\end{eqnarray*}
using Lemma \ref{Mlemma}. Therefore, we can apply the dominated convergence theorem to arrive at the limit,
\begin{equation}\label{Z}
Z_{t}^{0}(y)= \int_{0}^{t}\int_{U}P_{t-s}G(a,y,u_{s}^{0}(y))W(dads)
\end{equation}
The tightness result obtained above, implies that $\{Z_{t}^{\epsilon}\}$ is relatively compact by Prohorov's theorem. We apply the following tightness criterion stated in \cite{Whitt}.
\begin{theorem}[Corollary 11.6.1 in \cite{Whitt}]
Let $\{P_{n}\}_{n\geq 1}$ be a sequence of probability measures on a metric space $(S,m)$. If the sequence $\{P_{n}\}$ is tight and the limit of any convergence subsequence from $\{P_{n}\}$ must be $P$, then $P_{n}\xrightarrow{d} P$.
\end{theorem}
But every subsequence in our case has the form given in (\ref{mild}) and the uniqueness of solutions to (\ref{Z}) can be derived analogous to estimates in (\ref{unique}). Since $u_{s}^{0}(y)$ is a PDE then the integrand in Ito integral (\ref{Z}) is deterministic and by applying the H$\ddot{o}$lder inequality and following similar steps to (\ref{Zt}) we can show that
\begin{equation}\label{estimate}
\int_{0}^{1}\int_{U}\left(\int_{\mathbb{R}} P_{t-s}G(a,y,u_{s}^{0}(y))f(y)dy\right)^{2}\lambda(da)ds<\infty
\end{equation}
which implies that $Z_{t}^{0}$ is a Gaussian process with zero mean and covariance,
\begin{eqnarray}\label{Gaussian}
&&Cov(\left<Z_{t}^{0},f\right>,\left<Z_{t}^{0},g\right>\\
&=& \int_{0}^{t}\int_{U}\int_{\mathbb{R}}P_{t-s}G(a,y,u_{s}^{0}(y))f(y)dy\int_{\mathbb{R}}P_{t-s}G(a,r,u_{s}^{0}(r))g(r)dr
\lambda(da)ds\nonumber
\end{eqnarray}
 Note that because of the transition invariant property of the Lebesgue measure, (\ref{Gaussian}) is equivalent to (\ref{co1}).

\section{CLT for SBM and FVP}

Now we turn our attention to the two population models mentioned in the introduction. As for SBM, $\{\mu_{t}^{\epsilon}\}$, based on SPDE (\ref{SBM}), we have $G(a,y,u)= 1_{0\leq a \leq u}+1_{u\leq a \leq 0}$ which satisfies both conditions (\ref{con1}) and (\ref{con2}). Therefore, the tightness result obtained for the general SPDE can be used in this case. As for the limit, recall that $u_{t}^{\epsilon}(y)= \int_{0}^{y}\mu_{t}^{\epsilon}(dx)$ and thus for $f\in \mathcal{S}(\mathbb{R})$,
\begin{equation*}
<\mu_{t}^{\epsilon},f>= -<u_{t}^{\epsilon},f'>
\end{equation*}
then the centered functional process for SBM is found by
\begin{eqnarray}\label{SBMcenter}
\left<\tilde{Z}_{t}^{\epsilon}, f\right> &=& \left<\frac{1}{\sqrt{\epsilon}}\left(\mu_{t}^{\epsilon}-\mu_{t}^{0}\right), f\right>\\
&=& \left<\frac{1}{\sqrt{\epsilon}} \left(-u_{t}^{\epsilon}+u_{t}^{0}\right),f'\right>\nonumber\\
&=& -\left<Z_{t}^{\epsilon},f'\right>\nonumber\\
&=& -\int_{0}^{t}\int_{U}\int_{\mathbb{R}} P_{t-s}G(a,y,u_{s}^{\epsilon}(y))f'(y)dy W(dads)\nonumber
\end{eqnarray}
based on our estimates in the case of the general SPDE, the limit of (\ref{SBMcenter}) is
\begin{equation*}
\left<\tilde{Z}_{t}^{0}, f\right> = -\int_{0}^{t}\int_{U}\int_{\mathbb{R}} P_{t-s} G(a,y,u_{s}^{0}(y))f'(y)dyW(dads)
\end{equation*}
Namely, we have, 
\begin{eqnarray*}
&&\left<\tilde{Z}_{t}^{0},f\right>\\
&=&-\int_{0}^{t}\int_{\mathbb{R}}\int_{\mathbb{R}}\left(1_{0\leq a\leq u_{s}^{0}(y)}+1_{u_{s}^{0}(y)\leq a\leq 0}\right)f'(y)dyW(dsda)\\
&=& -\int_{0}^{t}\int_{0}^{u_{s}^{0}(y)}\int_{0}^{\infty}f'(y)dyW(dads) - \int_{0}^{t}\int_{u_{s}^{0}(y)}^{0}\int_{-\infty}^{0}f'(y)dyW(dads)\\
&=& -\int_{0}^{t}\int_{0}^{\infty}\int_{\left(u_{s}^{0}\right)^{-1}(a)}^{\infty}f'(y)dyW(dads) -\int_{0}^{t}\int_{-\infty}^{0}\int_{-\infty}^{\left(u_{s}^{0}\right)^{-1}(a)} f'(y)dyW(dads)\\
&=& \int_{0}^{t}\int_{0}^{\infty}f\left((u_{s}^{0})^{-1}(a)\right)W(dads)- \int_{0}^{t}\int_{-\infty}^{0}f\left((u_{s}^{0})^{-1}(a)\right)W(dads)\\
&=:& I_{t}^{1}(f)-I_{t}^{2}(f)
\end{eqnarray*}

then,
\begin{equation*}
Cov\left(<\tilde{Z}_{t}^{0}, f>,<\tilde{Z}_{t}^{0},g>\right)
= \mathbb{E}\left(\left(I_{1}^{t}(f)-I_{2}^{t}(f)\right)\left(I_{1}^{t}(g)-I_{2}^{t}(g)\right)\right)
\end{equation*}
Because of the measurability of $f$ and $g$ we can write $f\left((u_{s}^{0})^{-1}(a)\right)$ and $g\left((u_{s}^{0})^{-1}(a)\right)$ as limits of simple functions, $\sum_{i}\alpha_{i}1_{A_{i}}$ and $\sum_{j}\beta_{j}1_{B_{j}}$ respectively. Note that
\begin{eqnarray*}
&&\mathbb{E}\left(\int_{0}^{t}\int_{0}^{\infty}\sum_{i}\alpha_{i}1_{A_{i}}W(dsda)\right)\left(\int_{0}^{t}\int_{-\infty}^{0} \sum_{j}\beta_{j}1_{B_{j}}W(dads)\right)\\
&=& \sum_{i}\alpha_{i}\mathbb{E}\left(W\left(A_{i}\cap [0,t]\times [0,\infty)\right)\right)\sum_{j}\beta_{j}\mathbb{E}\left(W\left(B_{j}\cap [0,t]\times (-\infty,0]\right)\right)
\end{eqnarray*}
where we have used the independent scattered property of Gaussian measures in the last step. This yields to
\begin{equation*}
\mathbb{E}\left(I_{t}^{2}(f)I_{t}^{1}(g)\right)= \mathbb{E}\left(I_{t}^{1}(f)I_{t}^{2}(g)\right)=0
\end{equation*}

Moreover,
\begin{eqnarray*}
&&\mathbb{E}\left(I_{t}^{1}(f)I_{t}^{1}(g)\right)+ \mathbb{E}\left(I_{t}^{2}(f)I_{t}^{2}(g)\right)\\
&=& \int_{0}^{t}\int_{0}^{\infty}\left(f\left(\left(u_{s}^{0}\right)^{-1}(a)\right)g\left(\left(u_{s}^{0}\right)
^{-1}(a)\right)\right)dads\\
&&\hspace{.3cm}+ \int_{0}^{t}\int_{-\infty}^{0}\left(f\left(\left(u_{s}^{0}\right)^{-1}(a)\right)g\left(\left(u_{s}^{0}\right)
^{-1}(a)\right)\right)dads\\
&=& \int_{0}^{t}\left<\mu_{s}^{0},fg\right>ds
\end{eqnarray*}
where we used $y:=(u_{s}^{0})^{-1}(a)$ in the final step and observed the deterministic nature of the integrand. Notice that since $a\in \mathbb{R}$ then the integral with respect to $a$ has to be separated in this case; hence, we cannot use the properties of Ito integral to directly derive the covariance as we did for the general SPDE.

As mentioned in the introduction, Schied \cite{Schied} also achieved the CLT for SBM, however with a different setup. He considered the process
\begin{equation*}
\frac{1}{\beta}\left<f,X_{\beta^{2}t}-X_{0}P_{\beta^{2}t}\right>
\end{equation*}
where $f$ is a bounded Lipschitz continuous function in $\mathbb{R}^{d}$, $X_{t}$ is the SBM and comparing with our process, $\beta=\sqrt{\epsilon}$. After proving the tightness of this process in $\mathcal{C}\left([0,1];\mathbb{R}^{n}\right)$, he proved that the finite-dimensional marginal distributions converge weakly to those of an n-dimensional Wiener process $W$ with covariance,
\begin{equation}\label{Cov}
\mathbb{E}\left(W_{t}^{i}W_{t}^{j}\right)= 2t\int f_{i}f_{j}d\mu
\end{equation}
where $\mu$ is the initial measure of the SBM. This shows the weak convergence of the process to $W$ with covariance (\ref{Cov}). We note that the covariance of our limit is different since our process depends on $\epsilon$ as well as $t$.

As for FVP, from SPDE (\ref{FVP}) characterization, we can see that $G(a,y,u_{s}^{\epsilon}(y))= 1_{a\leq u_{s}^{\epsilon}(y)} - u_{s}^{\epsilon}(y)$, which also satisfies the two conditions (\ref{con1}) and (\ref{con2}) of $G(a,.,u_{s}^{\epsilon}(.))$ and so estimate (\ref{estimate}) holds. Thus we have,
\begin{equation*}
\left<\tilde{Z}_{t}^{0},f\right>= - \int_{0}^{t}\int_{0}^{1}\int_{\mathbb{R}}P_{t-s}\left(1_{a\leq u_{s}^{0}(y)}-u_{s}^{0}(y)\right)f'(y)dyW(dads)
\end{equation*}
is a Gaussian process with zero mean and covariance given below. For $f,g\in \mathcal{S}(\mathbb{R})$,
\begin{eqnarray*}
&&Cov\left(\left<\tilde{Z}_{t}^{0},f\right>,\left<\tilde{Z}_{t}^{0},g\right>\right)\\
&=& \int_{0}^{t}\int_{0}^{1}\int_{\mathbb{R}}\left(1_{a\leq u_{s}^{0}(y)}-u_{s}^{0}\right)f'(y)dy\int_{\mathbb{R}}\left(1_{a\leq u_{s}^{0}(r)}-u_{s}^{0}(r)\right)g'(r)drdads\\
&=& \int_{0}^{t}\int_{0}^{1}\left(\int_{(u_{s}^{0})^{-1}(a)}^{\infty}f'(y)dy + \left<\mu_{s}^{0},f\right>\right)\left(\int_{(u_{s}^{0})^{-1}(a)}^{\infty}g'(r)dr + \left<\mu_{s}^{0},g\right>\right)dads\\
&=& \int_{0}^{t}\int_{0}^{1}\left(-f\left((u_{s}^{0})^{-1}(a)\right)+ \left<\mu_{s}^{0},f\right> \right)\left(-g\left((u_{s}^{0})^{-1}(a)\right)+ \left<\mu_{s}^{0},g\right> \right)dads\\
&=& \int_{0}^{t}\int_{0}^{1}\left(-f(y)+\left<\mu_{s}^{0},f\right>\right)\left(-g(y)+ \left<\mu_{s}^{0},g\right>\right)\mu_{s}^{0}(dy)ds\\
&=& \int_{0}^{t}\int_{0}^{1}f(y)g(y)\mu_{s}^{0}(dy)ds - \int_{0}^{t}\int_{0}^{1}\left<\mu_{s}^{0},f\right>g(y)\mu_{s}^{0}(dy)ds\\
&& \hspace{.3cm} - \int_{0}^{t}\int_{0}^{1}f(y)\left<\mu_{s}^{0},g\right>\mu_{s}^{0}(dy)ds + \int_{0}^{t}\left<\mu_{s}^{0},f\right>\left<\mu_{s}^{0},g\right>ds
\end{eqnarray*}
Therefore, the CLT for FVP is achieved.

\section*{Acknowledgements}
I would like to thank my Ph.D. advisor, Dr. Jie Xiong for giving me useful ideas for this article. I am also thankful to Dr. Kei Kobayashi for pointing out the independent scattered property of Gaussian measures that was used in the case of SBM.

\end{document}